\documentclass{article}
\usepackage{amsthm,amstext,amsmath,amscd,amssymb,eucal,latexsym}
\usepackage{amscd}

\newcommand\F{\mathcal F}
\newcommand\E{\mathcal E}

\newcommand\NN{\mathbf N}

\newcommand\Ii{\mathcal I}

\newcommand\La{\mathcal L}

\newcommand\ZZ{\mathbf{Z}}

\newcommand\CC{\mathbf{C}}

\newcommand\Oc{\mathcal O}

\newcommand\Pq{\mathbf{P}^4}

\newcommand\Pj{\mathbf{P}}

\newtheorem*{RR}{Hirzebruch-Riemann-Roch Theorem}

\newtheorem{proposition}{Proposition}[section]
\newtheorem{definition}[proposition]{Definition}
\newtheorem{theorem}[proposition]{Theorem}
\newtheorem{question}[proposition]{Question}
\newtheorem{corollary}[proposition]{Corollary}
\newtheorem{remark}[proposition]{Remark}

\newtheorem{example}[proposition]{Example}

\newtheorem{lemma}[proposition]{Lemma}

\numberwithin{equation}{section}

\numberwithin{claim}{proposition}

\begin{document}

\title{Curves and vector bundles on quartic threefolds}
\author{Enrique Arrondo
\ and \ Carlo G. Madonna \footnote{The first author was supported
in part by the project BFM2003-03971/MATE funded by the Spanish
MCYT. The second author was supported in part by the projects I3P-CSIC and MTM2007-67623 founded by the Spanish
MEC.
}}
\date{}

\maketitle

\begin{abstract}
In this paper we study arithmetically Cohen-Macaulay (ACM for short) vector bundles $\E$
of rank $k \geq 3$ on hypersurfaces
$X_r \subset\Pj^4$ of degree $r \geq 1$.
We consider here
mainly the case of degree $r = 4$, which is the first
unknown case in literature. 
Under some natural conditions for
the bundle $\E$ 
we derive a list of possible Chern classes $(c_1,c_2,c_3)$
which may arise in the cases of rank $k=3$ and $k=4$, when $r=4$ and we give
several examples.
\end{abstract}

\section{Introduction}

Let $X=X_r \subset \Pj^4$ be a smooth hypersurface of
degree $r \geq 1$ and let $H$ be the class of a
hyperplane section. It is well known that
\begin{equation} \label{eq:ident}
{\rm Pic}(X_r) \cong \ZZ[H] \cong \ZZ
\end{equation}
and hence we may identify any line bundle $\mathcal L$ on
$X_r$ with the sheaf $\Oc_{X_r}(nH) \cong \Oc_{X_r}(n)$
for some $n \in \ZZ$. We recall the following:

\begin{definition}
Let $\E$ be a rank $k$ vector
bundle on $X_r$, $r \geq 1$. We call $\E$ arithmetically
Cohen-Macaulay (ACM for short) if
$h^i(X_r,\E(n))=0$
for all $i=1,2$ and $n \in \ZZ$, where $\E(n):=\E \otimes \Oc_{X_r}(n)$,
under the identification of (\ref{eq:ident}).
\end{definition}

In a previous work \cite{Ma1} the second author showed
that
a rank two ACM vector bundle
on $X_r$ splits as
a direct sum of two line bundles with only few possible exceptions.
Specifically, for each $r$, all the
possible first Chern classes of such an undecomposable
bundle (up to a twist with a line bundle) were given.
\par
For low values of $r$, i.e. $r=1,2,3,4$, a complete
classification of undecomposable rank two ACM vector
bundles is known. We refer the interested reader to
\cite{Hor} \cite{Ott} \cite{AC} \cite{Ma2} for more
details on these cases.
When $r=5$ the possible Chern classes of such an undecomposable
bundle where found in \cite{roum}
while, when $r \geq 6$ and $X_r$ is
general, the results of \cite{CM2} (for the case $r=6$) and of \cite{KRR} (for any $r \geq 6$) ensure the non
existence of such bundles.
Similar results were obtained for more general threefolds and
in particular a complete classification of rank two ACM
vector bundles on prime Fano threefolds of index $1$ was
given in
\cite{F} and \cite{AF}. Also in
some cases the corresponding moduli spaces are described (see
\cite{Ili} for a comprehensive account  of the known
results).
\par
On the other hand, for the higher rank case very little is
known even in the cases of low degree hypersurfaces. When
the rank is bigger than or equal to $3$, it is known that
there are no undecomposable ACM vector bundles if $r=1,2$,
while as far as we know, only the cubic case was considered
in
\cite{AC} under the following natural assumption (which we
generalize to arbitrary $r$):

\begin{definition}\label{def:star} We say that a rank
$k\geq 2$ vector bundle on $X_r \subset \Pj^4$ satisfies condition
$\star$ if the following conditions hold:\par
(i) $h^0\E(-1)=0$ and $h^0\E \geq k$;\par
(ii) there exist $k-1$ global sections of $\E$ whose
dependency locus has codimension equal to two;\par
(iii) $\E$ has not trivial summands.
\end{definition}

In this paper we will consider mainly the cases of low rank,
(precisely $k=3,4$) when $r=4$,
$X=X_4$ is general, and the above condition $\star$ is satisfied.
Notice, it is easy to produce examples of higher rank
($>3$) ACM vector bundles by using extension
classes. Specifically, starting with a pair $(\E',\E'')$ of
rank
$k$ and
$k'$ ACM bundles on $X_r$, one may always consider extensions
\begin{equation}\label{estensione}
0 \to \E' \to \E \to \E'' \to 0
\end{equation}
to get a rank $k'+k''$ ACM bundle on $X_r$.
Of course, since line bundles on $X_r$ are ACM, to get
such non trivial extension classes it is necessary to
have $k,k' \geq 2$.
We will give several examples of rank four ACM bundles
satisfying condition $\star$ when $r=4$ obtained in this way
in Section \ref{section:extensions}.
\par
Since such a construction was already used for the case
$r=3$ in \cite{AC} and for the case $r=5$ in \cite{MadCEJM},
one could state the following:

\begin{question} \label{q:1}
Let $\E$ be a rank four ACM vector bundle satisfying condition $\star$ on a general
hypersurface $X_r \subset \Pj^4$. Is $\E$
necessarily obtained as an extension class of a pair of
rank two ACM bundles on $X_r$, as in (\ref{estensione})?
\end{question}

We give here a negative answer to  the above question when
$r=4$ and $X_4$ is general (see Proposition \ref{ex:6-3}).
We notice that Question \ref{q:1} still remains open when
$r=3$.
Moreover when $r=6$ we have not examples at all of rank
$4$ ACM bundles, while when $r=5$ we know actually
only one example, given in \cite{MadCEJM}, which satisfies condition
$\star$.
\par\medskip
Turning to the case of rank $k=3$, we want to answer the following
question:
\begin{question} \label{q:2}
Are there rank three ACM bundles satisfying condition $\star$ on a general
hypersurface $X_r \subset \Pj^4$ of degree $r \geq 4$?
\end{question}

Since in \cite{AC} an affirmative answer (see Proposition 4.9) was already given for the case $r=3$ we
consider the case $r = 4$ here. We give an affirmative
answer
to this last question when $r=4$. We notice that Question
\ref{q:2} remains open when $r \geq 5$.
To answer the previous questions, we
derive the possible first Chern class of an ACM bundle of
rank $k \geq 3$ satisfying condition $\star$ on any smooth
hypersurface of degree $r \geq 3$. Then we consider in
detail the cases of low rank, i.e. $k=3$ and $k=4$, when
$r=4$, and we go on by a case by case analysis.
In this direction our main result is the following:

\begin{theorem} \label{thm:primoelenco}
Let $\E$ be an ACM vector bundle of rank $k
\in \{3,4\}$ on a general quartic threefold
$X_4\subset\Pj^4$ which satisfies condition $\star$.
Then the possible Chern classes $c_i=c_i(\E)$ and the genus
$g=g(C)$ of the curve $C$, dependency locus of $k-1$ global
section of $\E$, are given in the following table:
$$
\begin{matrix}
k & c_1 & c_2 \in & c_3 & g  \cr
\hline\vspace{-3.8mm}
\cr
3 & 1 & 5 & 2 & 2 \cr
3 & 2 & [8,11] & c_2-6 & c_2-2  \cr
3 & 3 & [17,18] & 2c_2-26 & 2c_2-12 \cr
3 & 4 & [27,28] & 3c_2-66 & 3c_2-32  \cr
4 & 1 & 6 & 4 & 3 \cr
4 & 2 & [8,12] & c_2-4 & c_2-1 \cr
4 & 3 & [16,22] & 2c_2-24 & 2c_2-11 \cr
4 & 4 & [28,32] & 3c_2-64 & 3c_2-31 \cr
4 & 5 & [44,46] & 4c_2-132 & 4c_2-65 \cr
4 & 6 & 64 & 84 & 203  \cr
\end{matrix}
$$
Moreover there are examples of bundles with invariants $(k;c_1,c_2,c_3)$
for any of the following quadruples:
$$
\begin{aligned}
{}& (4;1,6,4), (4;2,\alpha,\alpha-4),
(4;3,\beta,2\beta-24),
(4;4,\gamma,3\gamma-64),\\
& (4;5,46,52),
(4;6,64,84),(3;1,5,2), 
\end{aligned}
$$
where $\alpha \in \{10,11,12\}$, $\beta \in
\{19,20 \}$, and $\gamma \in \{29,30,32 \}$.
\end{theorem}

In the case of the rank four, almost all the examples are given
by considering extension classes of rank two ACM bundles.
It is not so for the case $(4;1,6,4)$, which is
constructed by means of the Hartshorne-Serre correspondence
(see Proposition \ref{prop:vogelaar})
starting from a suitable smooth projectively normal
curve in $X_4$ and it is not given by any extension class of any two rank
two ACM bundles on $X_4$.
The example for the rank three case is given in similar
way, by showing the existence of the corresponding curve
with prescribed invariants and then applying the above
mentioned Hartshorne-Serre correspondence.
\par
The paper is organized as follows. In Section \ref{section:facts}
we recall some known facts and generalities that we need in the
paper. In Section \ref{section:extensions} (see Lemma
\ref{cor:ext}), we find conditions such that the direct sum (or
more generally any extension class) of ACM vector bundles
satisfying condition $\star$ is still an ACM vector bundle
satisfying condition $\star$, and produce from this several
examples of rank four bundles on $X_4$ (precisely all but one in
the list of Theorem \ref{thm:primoelenco}). In Section
\ref{section:curves} we use the Hartshorne-Serre correspondence
(see Proposition \ref{prop:vogelaar}) between curves and vector
bundles to give a characterization (see Proposition
\ref{pro:criterio-ACM}) of when the vector bundle obtained from a
curve is ACM and satisfies condition $\star$; using this, we
construct the example of rank three and the remaining one of
rank four of ACM vector bundles satisfying condition $\star$ on
$X_4$ of the list of Theorem \ref{thm:primoelenco}. Finally, in
section \ref{sec:bounds} we give restrictions for the Chern
classes of ACM vector bundles on $X_r$ satisfying condition
$\star$, which when particularized for $r=4$ and rank $k=3,4$
yields the list given in Theorem \ref{thm:primoelenco} of all
possible Chern classes in these ranks.
\par
The preliminary draft of this paper was completed during the
second author stay in the Departamento de \'Algebra of University
Complutense in Madrid, supported by a grant ``Mensilit\`a per
l'estero" from the ``Istituto Nazionale di Alta Matematica
F.Severi" (10/2005-03/2006).

\section{Preliminaries and basic facts} \label{section:facts}
In this paper we work over the field of complex numbers $\CC$. Let
us consider a smooth hypersurface
$X_r \subset \Pj^4$
of degree
$r \geq 1$ in
the complex 4-dimensional projective space $\Pj^4$.
Then, as recalled in the introduction
${\rm Pic}(X_r) \cong \ZZ[H]$
where $H$ is the class of a hyperplane section,
and the canonical divisor of $X_r$ is
$K_{X_r}=(r-5)H$.
\par
If $P$ denotes the class of a point and $L$ the class of
a line, then the intersection products on $X_r$ are given by
$H^3=rP=r \ , \ HL=1P=1 \ , H^2=rL$.
Let $\E$ be a rank $k$ vector bundle on $X_r$. We
identify the first Chern class $c_1(\E)$ of $\E$ with the
integer number $c_1$ which corresponds to $c_1(\E)$ under
the above isomorphism, i.e.
$c_1(\E)=c_1H=c_1.$
In a similar way we identify
$c_2(\E)=\deg c_2(\E)=c_2$
and
$c_3(\E)=c_3P=c_3.$
Under this identification to any rank $k$ vector bundle $\E$ on
$X_r$ corresponds a quadruple
$(k;c_1,c_2,c_3) \in \NN \times \ZZ^3$.
In this language, our main result
gives all the possible quadruples of a rank three and four ACM vector
bundle on a smooth quartic threefold which could arise under the
condition $\star$ of Definition \ref{def:star}.
\par
For further computations we write down explicitly the
Chern classes of the bundle $\E(n)=\E \otimes \Oc_{X_r}(n)$ in the following
equations:
\[
\begin{aligned}
c_1(\E(n)) & =c_1+kn \\
c_2(\E(n)) &=c_2+r(k-1)nc_1+r\binom{k}{2}n^2=\\
& =c_2+rn(k-1)\big(c_1+\frac12 nk\big) \\
c_3(\E(n)) &=c_3+(k-2)nc_2+r\binom{k-1}{2} n^2c_1
+r\binom{k}{3} n^3=\\
& =c_3+(k-2)n\big(c_2+\frac12 (k-1)nr c_1+\frac16
rn^2k(k-1)\big)
\end{aligned}
\]
In the sequel, to perform some computations, we will
frequently use the following version of Riemann-Roch
theorem for vector bundles:

\begin{RR} If $\E$ is a rank $k$ vector bundle on a
smooth hypersurface $X_r \subset \Pj^4$ of degree $r \geq
1$ with Chern classes $c_i(\E)=c_i \in \ZZ$ for $i=1,2,3$,
then
\begin{equation} \label{eq:RRk}
\begin{aligned}
\chi(\E) &=\frac16 rc_1^3-\frac12 c_1c_2+\frac12
c_3+\frac14 rc_1^2(5-r) -\frac12 c_2(5-r)+\\
& +
\frac{1}{12}\big((r-5)^2+(10-5r+r^2)\big)c_1r+\frac1{24}rk(5-r)(10-5r+r^2)
\end{aligned}
\end{equation}
In particular 
\begin{equation} \label{eq:RR1}
\begin{aligned}
\chi(\Oc_{X_r}(a)) &=\frac16 a^3r+\frac14 a^2 r(5-r) +
\frac1{12}ar \big((r-5)^2+ \\ 
& + (10-5r+r^2)\big)
+ \frac{r}{24}(5-r)(10-5r+r^2)
\end{aligned}
\end{equation}
\end{RR}
Recall that a vector bundle $\E$ is ``decomposable'' if
there exist $\E'$ and $\E''$ such that
$\E \cong \E' \oplus \E''.$
Otherwise $\E$ is ``undecomposable''.
By the Serre's duality, $\E$ is
ACM if and only if
\[
h^1\E(n)= h^1 \E^{\vee}(-n+r-5)=0 \ \ \forall n \in \ZZ
\]
which is equivalent to the condition
$h^1\E(n)=h^1\E^{\vee}(n)=0 \ \ \forall n \in \ZZ$.
\par
From now on we will consider mainly vector bundles
$\E$ of rank $k \geq 3$ for which the condition $\star$
(see Definition \ref{def:star}) is
satisfied.
If $h^0\E(-1)=0$ and $h^0\E>0$ we will say that $\E$ is
``normalized''. It is always possible to assume this
condition since we may replace $\E$ with $\E(-b)$, where
$b=b(\E)$ is defined as  \begin{equation}\label{def-b}
b(\E)=b=\max \{n \in \ZZ: h^0\E(-n)=0 \}.
\end{equation}
\noindent We recall also the notion of stability that we will use in the sequel. A rank
$k$ vector bundle $\E$ on $X_r$ is ``stable'' (resp.
``semistable'') if for any  subbundle $\F \subset\E$ of rank
$m<k$ we have
$c_1(\F)/m < c_1(\E)/k$
(resp. $\leq$). In particular if $b(\E)=b$ is defined as
in (\ref{def-b}) we have $kb(\E)-c_1(\E)<0$ if $\E$ is
stable (resp. $\leq$). Indeed, any global section of
$H^0\E(-b)$ gives rise to a morphism $\Oc_{X_r} \to
\E(-b)$ and hence if $\E$ is stable then
$0<c_1(\E(-b))/k=(c_1(\E)-kb)/k$ and hence $kb-c_1<0$. In
particular if $\E$ is normalized and stable then $c_1>0$.
When $k=2$ the condition $kb(\E)-c_1(\E) \leq 0$ is in
fact equivalent to the semistablity of $\E$. If $k \geq 3$ this is not an equivalence any more.
\par
We can use condition (ii) in Definition \ref{def:star} to
translate properties of vector bundles satisfying $\star$ to
properties of curves in $X_r$. Indeed, for a vector bundle
$\E$ satisfying condition $\star$, the choice of $k-1$ global
sections whose dependency locus is a curve
$C\subset X_r$ yields the fundamental exact sequence:
\begin{equation} \label{eq:basic}
0 \to \Oc_{X_r}^{k-1} \to \E \to \Ii_{C}(c_1) \to 0
\end{equation}
where $\Ii_C$ is the ideal sheaf of $C \subset X_r$. When
dualizing (\ref{eq:basic}), we get the exact sequence
\begin{equation} \label{eq:basic3}
0 \to \Oc_{X_r}(-c_1) \to \E^{\vee} \to \Oc_{X_r}^{k-1}
\to \omega_C(5-r-c_1) \to 0
\end{equation}
where $\omega_C\cong{\cal E}xt^2(\Oc_C,\omega_{X_r})$ is the
dualizing sheaf of $C$. This implies, in particular, that
$\omega_C(5-r-c_1)$ is generated by $k-1$ global sections. Next we
recall the generalized Hartshorne-Serre correspondence (see
\cite{vog} \cite{arrondo}), which states that is possible to
reverse this process, in the sense that one recovers the vector
bundle from the surjection in (\ref{eq:basic3}):

\begin{proposition}\label{prop:vogelaar}
Let $C\subset X_r$  be a locally complete intersection curve
and assume that, for some $c_1\in\ZZ$,
$
\omega_C(5-r-c_1)
$
has $k-1$ generating
global sections. Then, there exists a rank $k$ vector bundle
$\E$
on $X_r$ fitting in the exact sequence (\ref{eq:basic}) and the
surjection of (\ref{eq:basic3}) corresponds to the choice of the $k-1$
global sections of $\omega_C(5-r-c_1)$. In particular
$c_1(\E)=c_1$,
$c_2(\E)=\deg C$
and $C$ is the dependency locus of the $k-1$ global sections
of $\E$ given by (\ref{eq:basic}).
\end{proposition}

This gives, any time we have an exact sequence like
(\ref{eq:basic}), the relation between the two first Chern classes
of $\E$ and the degree of $C$. Moreover, the third Chern
class of $\E$ is determined by the (arithmetic) genus of the
curve $C$:

\begin{corollary}
If $C \subset X_r$ is a locally complete intersection curve
associated to a vector bundle on $X_r$ of rank $k
\geq 2$ as above, then
\[
g(C)=-\frac52 c_2+\frac12 c_1c_2+\frac12
c_3+\frac{25}{12} r+\frac12 rc_2-\frac{35}{24} r^2
+\frac{5}{12} r^3-\frac1{24} r^4
\]
where $c_i=c_i(\E)$, $i=1,2,3$. In particular when $r=4$ we have
\begin{equation} \label{eq:genus}
g(C)=1+\frac12 c_1c_2-\frac12 c_2+\frac12 c_3.
\end{equation}
\end{corollary}

\begin{proof}
It follows from the equalities
$g(C)=1-\chi(\Oc_C)=1-\chi(\Oc_{X_r})+\chi(\Ii_C)$,
$\chi(\Ii_C)=\chi(\E(-c_1))-\chi(\Oc_{X_r}^{k-1}(-c_1))$
(by the exact sequence (\ref{eq:basic})) and the
Riemann--Roch theorem.
\end{proof}

Having in mind Hartshorne-Serre correspondence, we will
frequently consider pairs
\begin{equation} \label{eq:pair}
(\E,C)
\end{equation}
given by a vector bundle $\E$ of rank $k \geq 3$ which
satisfies condition $\star$ and a curve $C
\subset X_r$ dependency locus of $k-1$ global sections of
$\E$. We will also say in this case that the pair $(\E,C)$
is a {\it ACM pair} and satisfies condition $\star$ (in Proposition
\ref{pro:criterio-ACM} we will give a criterion for this
property in terms of the curve $C$).

We recall the following generalization of Clifford's theorem
(observe that in our situation $C$ is connected because the
fact that $\E$ is ACM implies from (\ref{eq:basic}) the
vanishing of $H^1 \Ii_C$):

\begin{theorem} \label{theo:Clifford} Let $C$ be a
connected locally complete intersection curve, and let $\La$
be a line bundle on $C$ such that $H^0 \La \ne0$ and $H^1 \La \ne0$.
Then $2(h^0 \La -1)\le\deg \La$.
\end{theorem}

\begin{proof} The same proof as in \cite{cl} Theorem A works
in this case, since the irreducibility assumed there is not used
to prove the inequality. In fact, we still have a nondegenerate
bilinear map $H^0 \La \otimes Hom(\La,\omega_C)\to H^0 \omega_C $, so
that the result follows, as in  \cite{cl}, by the so-called
bilinear lemma, Serre's duality and Riemann-Roch theorem. The
connectedness of $C$ is used to conclude that the dimension
$h^0 \omega_C=h^1 \Oc_C$ is $g(C)$ because $h^0 \Oc_C=1$.
\end{proof}

Let us finish this section by recalling the classification of
undecomposable rank two ACM bundles on quartics 
$X_4$ obtained in \cite{Ma2} and on cubics $X_3$
obtained in \cite{AC}.

\begin{proposition} \label{prop:list}
Let $\E$ be a normalized and undecomposable rank two ACM
bundle on a smooth hypersurface $X_r \subset \Pj^4$ with Chern classes $c_i(\E)=c_i$,
$i=1,2$.
\begin{enumerate}
\item If $r=4$ then
$(c_1,c_2) \in \{
(-1,1),(0,2),(1,3),(1,4),(1,5),(2,8),(3,14) \}$;
\item if $r=3$ then
$(c_1,c_2)=\{(0,1),(1,2),(2,5) \}$.
\end{enumerate}
Moreover all the cases arise.
\end{proposition}

\begin{remark} \label{rem:lemma}
\normalfont{
In particular in all the above cases condition $\star$ is satisfied with the only exceptions of
cases $(c_1,c_2) \in \{(-1,1),(0,2),(1,5) \}$ when $r=4$ and
$(c_1,c_2)=(0,1)$ when $r=3$. Moreover, as shown in
\cite{Ma2}, it also holds that the vector bundles on
$X_4$ with $(c_1,c_2)=(3,14)$  and the general one with
$(c_1,c_2)=(2,8)$ are generated by their global sections.
}\end{remark}

\section{Constructing bundles from extensions}
\label{section:extensions}

We prove in this section that direct sums (and more
generally extensions) of vector bundles satisfying
condition $\star$ are still ACM vector bundles
satisfying condition $\star$. We will end by producing in
this way examples of rank four ACM vector bundles satisfying
condition $\star$ on $X_4$.
Observe first that an ACM vector bundle $\E$ satisfying
condition $\star$ cannot decompose as
$\E \cong \E' \oplus \Oc_{X_r}(a)$
where $\E'$ is a vector bundle of rank $k-1$. Indeed, the
possibility $a>0$ is excluded by condition (i), the
possibility $a=0$ is excluded by condition (iii) and  the
possibility $a<0$ is excluded by condition (ii), since the
dependency locus of $k-1$ sections of $\E$ would be the
dependency locus of $k-1$ sections of $\E'$ (hence of
expected codimension one).

Of course it can be possible to have direct sums
$\E \cong \E' \oplus \E''$ with $rk(\E')$ and  $rk(\E'') \geq 2$.
In fact, we are going to see that the second part in
condition (i) implies that condition $\star$ is preserved by
direct sums (and more generally by extensions). To show
this, we start with a standard fact:

\begin{lemma} \label{lemma:sum}
Let $\E'$ and $\E''$ be two vector bundles on $X_r$ of rank $k'$
and $k''$ respectively. Let $s'_1,\ldots,s'_{k'}$ be linearly
independent sections of $\E'$ and let $s''_1,\ldots,s''_{k''}$ be
linearly independent sections of $\E''$. Assume that the
dependency locus of $s'_1,\ldots,s'_{k'-1}$ is a curve $C'$ and
the dependency locus of $s''_2,\ldots,s''_{k''}$ is a curve $C''$.
Then the dependency locus of the sections
\[
(s'_1,0),\ldots,
(s'_{k'-1},0),(s'_{k'},s''_1),(0,s''_2),\ldots,(0,s''_{k''})
\]
of
$\E'\oplus\E''$ is $C'\cup
C''\cup V(s'_1\wedge\ldots\wedge
s'_{k'},s''_1\wedge\ldots\wedge s''_{k''})$, where
$s'_1\wedge\ldots\wedge s'_{k'}$ is interpreted as a section
of the line bundle $\wedge^{k'}\E'$ and
$s''_1\wedge\ldots\wedge s''_{k''}$ as a section of
$\wedge^{k''}\E''$.
\end{lemma}

\begin{proof}
It is enough to check the statement locally. We can thus
restrict ourselves to an open subset of $X_r$ on which the
vector bundles trivialize. If on that open set the section
$s'_i$ is represented by the functions
$f'_{i1},\ldots,f'_{ik'}$ and the section $s''_j$ is
represented by
$f''_{j1},\ldots,f''_{jk''}$, then the dependency locus of
$(s'_1,0),\ldots,
(s'_{k'-1},0),(s'_{k'},s''_1),(0,s''_2),\ldots,(0,s''_{k''})$
is defined by the vanishing of the maximal minors of the
$(k'+k''-1)\times(k'+k'')$matrix
$$
\begin{pmatrix}
f'_{11}&\ldots&f'_{1k'}&0&\ldots&0\cr
\vdots&&\vdots&\vdots&&\vdots\cr
f'_{k'-1,1}&\ldots&f'_{k'-1,k'}&0&\ldots&0\cr
f'_{k'1}&\ldots&f'_{k'k'}&f''_{11}&\ldots&f''_{1k''}\cr
0&\ldots&0&f''_{21}&\ldots&f''_{2k''}\cr
\vdots&&\vdots&\vdots&&\vdots\cr
0&\ldots&0&f''_{k''1}&\ldots&f''_{k''k''}
\end{pmatrix}
$$
These $k'+k''-1$ minors take the form
$$
\left|\begin{matrix}
f'_{11}&\ldots&f'_{1,i-1}&f'_{1,i+1}&\ldots&f'_{1k'}\cr
\vdots&&\vdots&\vdots&&\vdots\cr
f'_{k'-1,1}&\ldots&f'_{k'-1,i-1}&f_{k'-1,i+1}&\ldots&f'_{k'-1,k'}
\end{matrix}\right|
\left|\begin{matrix}
f''_{11}&\ldots&f''_{1k''}\cr
\vdots&&\vdots\cr
f''_{k''1}&\ldots&f''_{k''k''}
\end{matrix}\right|,
\ \ \ \ \ \ \ \ \ \ \ i=1,\ldots,k'
$$
$$
\left|\begin{matrix}
f'_{11}&\ldots&f'_{1k'}\cr
\vdots&&\vdots\cr
f'_{k'1}&\ldots&f'_{k'k'}
\end{matrix}\right|
\left|\begin{matrix}
f''_{21}&\ldots&f''_{2,j-1}&f''_{2,j+1}&\ldots&f''_{2k''}\cr
\vdots&&\vdots&\vdots&&\vdots\cr
f''_{k''1}&\ldots&f''_{k'',j-1}&f''_{k'',j+1}&\ldots&f_{k''k''}
\end{matrix}\right|,
\ \ \ \ \ \ \ \ \ \ \ j=1,\ldots,k''
$$
Since the curve $C'$ is locally defined by the minors
$$
\left|\begin{matrix}
f'_{11}&\ldots&f'_{1,i-1}&f'_{1,i+1}&\ldots&f'_{1k'}\cr
\vdots&&\vdots&\vdots&&\vdots\cr
f'_{k'1}&\ldots&f'_{k',i-1}&f'_{k',i+1}&\ldots&f'_{k'k'}
\end{matrix}\right|,\ \ \ \ \ \ \ \ \ i=1,\ldots,k',
$$
the curve $C''$ is
defined by the minors
$$
\left|\begin{matrix}
f''_{11}&\ldots&f''_{1,j-1}&f''_{1,j+1}&\ldots&f''_{1k''}\cr
\vdots&&\vdots&\vdots&&\vdots\cr
f''_{k'1}&\ldots&f''_{k'',j-1}&f''_{k'',j+1}&\ldots&f''_{k''k''}
\end{matrix}\right|,\ \ \ \ \ \ \ \ j=1,\ldots,k''
$$
and the set
$$
V(s'_1\wedge\ldots\wedge s'_{k'},
s''_1\wedge\ldots\wedge s''_{k''})
$$
is defined by
$$
\left|\begin{matrix} f'_{11}&\ldots&f'_{1k'}\cr
\vdots&&\vdots\cr
f'_{k'1}&\ldots&f'_{k'k'}
\end{matrix}\right|
$$
and
$$
\left|\begin{matrix}
f''_{11}&\ldots&f''_{1k''}\cr
\vdots&&\vdots\cr
f''_{k''1}&\ldots&f''_{k''k''}
\end{matrix}\right|,
$$
the result follows at once.
\end{proof}

As a corollary, we can prove the following (see also Lemma
4.1 of \cite{AC}):

\begin{lemma} \label{lemma:ext}
Let $\E'$ and $\E''$ be two vector bundles on $X_r$ of rank
$k'$ and $k''$ respectively. Assume that $\E'$ and $\E''$
satisfy condition $\star$ and that there are sections
$s'_1,\ldots,s'_{k'}$ of $\E'$ and $s''_1,\ldots,s''_{k''}$
of
$\E''$ such that the hypersurfaces
$V(s'_1\wedge\ldots\wedge s'_{k'})$ and
$V(s''_1\wedge\ldots\wedge s''_{k''})$ do not share a common
component. Then a general $\E$ fitting in an exact sequence
\begin{equation}\label{eq:extension}
0\to\E'\to\E\to \E''\to0
\end{equation}
satisfies condition $\star$.
\end{lemma}

\begin{proof}
It is clear that $\E$ satisfies conditions (i) and
(iii) of Definition \ref{def:star}. Hence it is enough to
check condition (ii). Since this condition is open, it
suffices to prove it for $\E=\E'\oplus\E''$. But this follows
immediately from Lemma \ref{lemma:sum}, since our assumption
implies that we can find sections $s'_1,\ldots,s'_{k'}$ of
$\E'$ and $s''_1,\ldots,s''_{k''}$ of $\E''$ such that
$V(s'_1\wedge\ldots\wedge s'_{k'},s''_1\wedge\ldots\wedge
s''_{k''})$ is a curve.
\end{proof}

\begin{remark} \normalfont{We can use condition (ii) in
Definition \ref{def:star} to give a geometric interpretation
(and a criterion) for the hypothesis in Lemma
\ref{lemma:ext}. For a vector bundle $\E$ satisfying
condition $\star$, the choice of $k-1$ global sections
$s_1,\ldots,s_{k-1}$ whose dependency locus is a curve
$C\subset X_r$ yields the exact sequence
(\ref{eq:basic}). Since we are assuming $h^0 \E \ge k$, this
means that there is at least a section $s_k$ of $\E$
independent of $s_1,\ldots,s_{k-1}$. Any such $s_k$ yields a
hypersurface of degree $c_1$ containing $C$ and defined by
the global section
$s_1\wedge\ldots\wedge s_k$ of $\Oc_{X_r}(c_1)$ (with the
convention of Lemma \ref{lemma:sum}).
Notice if $\E'$ and $\E''$ are ACM also $\E$ is ACM. 
}\end{remark}

\begin{corollary} \label{cor:ext}
Let $(\E',C')$ and $(\E'',C'')$ be two
ACM pairs satisfying condition $\star$. Then a general $\E$
appearing in an extension as (\ref{eq:extension}) satisfies
condition $\star$ if at least one of the following
conditions is satisfied:
\begin{enumerate}
\item[(a)] there is a hypersurface of degree $c'_1$
containing $C'$ and a hypersurface of degree $c''_1$ containing $C''$ such that these two hypersurfaces do not share any
component;
\item[(b)] at least one of $\E'$ and $\E''$ is generated by
its global sections.
\end{enumerate}
\end{corollary}

\begin{proof}
Part (a) is an immediate consequence of Lemma
\ref{lemma:ext} and the previous remark. To prove part
(b), assume for instance that $\E'$ is generated by its
global sections. We fix any hypersurface of degree $k''$
containing $C''$ and take a point in any of its components.
Since $\E'$ is generated by its global sections we can find
$k'$ sections of it whose dependency locus does not contain
any of those points. This dependency locus is therefore a
hypersurface of degree $c'_1$ containing $C'$ and having no
common components with the fixed hypersurface of degree
$c''_1$ containing $C''$. We thus conclude from (a).
\end{proof}

We want now to apply Corollary \ref{cor:ext} to obtain rank
four ACM vector bundles on $X_4$ satisfying condition $\star$
from vector bundles of rank two $\E'$ and $\E'$. We recall
how to compute the invariants of any $\E$ fitting in an
extension (\ref{eq:extension}). Let us denote by $c'_i \in
\ZZ$ the Chern classes of $\E'$ and by $c''_i \in \ZZ$ the
Chern classes of $\E''$. Then:
\begin{equation} \label{eq:wit}
\begin{aligned}
c_1(\E) & =c'_1+c''_1 \\
c_2(\E) & =c'_2+4 c'_1 c''_1 +c''_2 \\
c_3(\E) & =c'_2c''_1+c'_1c''_2
\end{aligned}
\end{equation}
Checking Remark \ref{rem:lemma} for finding the list of rank two ACM bundles
on $X_4$ satisfying condition $\star$ we find the following
list of examples.

\begin{example}
\label{ex:4,6,64} \normalfont{ We take $\E',\E''$ to be rank two
ACM vector bundles on $X_4$ with $c_1=3$, $c_2=14$. These are
generated by their global sections, so that it follows from
Corollary \ref{cor:ext}(b) that any general element in
$Ext^1(\E'',\E')$ (which is a vector space of dimension at least
$7$) provides a rank four ACM vector bundle satisfying condition
$\star$. In fact, in this case it is easier to observe that any
extension provides a globally generated vector bundle, and hence
it always satisfies condition $\star$. By the formulas
(\ref{eq:wit}), the invariants of such vector bundle are
$(k;c_1,c_2,c_3)=(4;6,64,84).$
}\end{example}

\begin{example}
\label{ex:4,5,46} \normalfont{ We take now the rank two ACM vector
bundles $\E',\E''$ in $X_4$ with $c_1(\E')=3$, $c_2(\E')=14$, and
$c_1(\E'')=2$. It follows again from Corollary \ref{cor:ext}(b)
that a general element in $Ext^1(\E'',\E')$ (which is a vector
space of dimension at least $7$) provides a rank four ACM vector
bundle $\E$ satisfying condition $\star$. In this case, the
invariants produced by the formulas (\ref{eq:wit}) are
$(k;c_1,c_2,c_3)=(4;5,46,52).$
}\end{example}

\begin{example} \label{ex:4,4,293031}
\normalfont{
We repeat the same reasoning as in Example \ref{ex:4,5,46},
but taking now $\E''$ a rank two ACM vector bundle with
$c_1(\E'')=1$, $c_2(\E'')=d$ (for $d=3,4$). As before, a
general element in $Ext^1(\E'',\E')$ (which is a vector
space of dimension at least $3d-8$) yields an unstable rank
four ACM vector bundle satisfying condition $\star$ and with
invariants
$(k;c_1,c_2,c_3)=(4;4,26+d,3d+14), \, \, d=3,4.$
}\end{example}

\begin{example}\label{ex:4,4}
\normalfont{
Now we take  $\E',\E''$ two normalized rank two ACM vector
bundles that are globally generated on $X_4$ with $c_1=2$ and
$c_2=8$. Again from Corollary \ref{cor:ext}(b) we get that a
general element in $Ext^1(\E'',\E')$ (which is a vector space
of dimension at least $4$) provides a rank four ACM vector
bundle satisfying condition $\star$. Using once more
the equations (\ref{eq:wit}) we find that this bundle has
invariants
$(k;c_1,c_2,c_3)=(4;4,32,32).$
}\end{example}

\begin{example} \label{ex:4,3}
\normalfont{
We replace now in Example \ref{ex:4,4} the vector bundle
$\E''$ with a normalized rank two ACM vector bundle on $X_4$ with
$c_1(\E'')=1$ and $c_2(\E'')=d$ (with $d=3,4$). Hence we get that
a general element in $Ext^1(\E'',\E')$ (which is a vector space of
dimension at least $2d-4$) provides an unstable rank four ACM
vector bundle satisfying condition $\star$ and with invariants
$(k;c_1,c_2,c_3)=(4;3,16+d,8+2d), \, \, d=3,4.$
}\end{example}

\begin{example} \label{ex:4,2}
\normalfont{
Take now $\E'$ and $\E''$ to be two normalized rank
two ACM bundles on $X_4$ with $c_1(\E')=1$,
$c_2(\E')=d'$ and $c_1(\E'')=1$, $c_2(\E'')=d''$,
(with $d', d''=3,4$). The curves obtained as the
zero loci of sections of $\E'$ and $\E''$ are
elliptic curves of degree $d'$ and $d''$, in any
case contained in at least one hyperplane of
$\Pj^4$. A simple calculation shows that a
hyperplane section of $X_4$ (which is a K3
surface) contains at most a pencil of elliptic
curves of degree $d$, while the family of those
curves in $X_4$ has dimension $d$ (see \cite{Ma2}).
Hence we  can take pairs $(\E',C')$ and $(\E'',C'')$
such that the hyperplanes containing $C'$ and $C''$
are different. Therefore we can apply Lemma
\ref{lemma:ext} and take a general element in
$Ext^1(\E',\E'')$ (which has dimension at least
$d'+d''-6\ge0$) to produce an ACM vector bundle
satisfying condition $\star$ and with invariants
$(k;c_1,c_2,c_3)=(4;2,d'+d''+4,d'+d''), \, \, d,d'=3,4.$
}\end{example}

\section{Constructing bundles from curves on the quartic
threefold}
\label{section:curves}

In this section, we use first Hartshorne-Serre
correspondence, to translate the property of being an ACM
vector bundle satisfying condition $\star$ to the property of curves in
$X_r$ to be associated to it. This will allow to produce several examples of ACM
vector bundles of rank $k=3,4$ on $X_4$ satisfying condition
$\star$ from curves in $X_4$.

We start by characterizing when a curve $C$ determines an ACM
vector bundle (observe that the condition $c_1>0$ will not be
restrictive because of (\ref{primobound:c1})). In the sequel
we will denote by $\Oc_C(1)$ the restriction of the
hyperplane class $H$ of $X_r$ to $C \subset X_r$.

\begin{proposition} \label{pro:criterio-ACM}
Let $C\subset X_r$ be a locally complete intersection curve
and let $\E$ be the vector bundle of rank $k \geq 2$
obtained, as in Proposition \ref{prop:vogelaar}, from $k-1$
sections of $\omega_C(5-r-c_1)$. Assume $c_1>0$. Then $\E$ is
an ACM bundle satisfying condition $\star$ if and only if
$C$ is projectively normal and the following four conditions
hold:
\begin{enumerate}
\item[(a)] the $k-1$ given sections form a basis of
$H^0\omega_C(5-r-c_1)$ (in particular
$h^0\omega_C(5-r-c_1)=k-1$);
\item[(b)] $h^0\Ii_C(c_1-1)=0$ and $h^0\Ii_C(c_1)\ge1$;
\item[(c)] $h^0\omega_C(4-r-c_1)=0$, for which a sufficient condition is
$2g(C)-2<(r+c_1-4)\deg C$;
\item[(d)] the natural map $\alpha_n:H^0 \omega_C(5-r-c_1) \otimes
H^0 \Oc_C(n)\to H^0\omega_C(5-r-c_1+n)$ is surjective
$\forall n \geq 1$.
\end{enumerate}
\end{proposition}

\begin{proof}
Since $h^i\Oc_{X_r}(n)=0$ for all $n \in\ZZ$ and $i=1,2$ then
the condition $h^1\E(n)=0$ for all $n \in\ZZ$ is equivalent
to the condition $h^1\Ii_C(c_1+n)=0$ for all $n \in\ZZ$, i.e.
to the projective normality of $C$. On the other hand,
splitting (\ref{eq:basic3}) into the following exact
sequences
\begin{equation}\label{eq:dual-sinistra}
0 \to \Oc_{X_r}(-c_1) \to \E^{\vee} \to K \to 0
\end{equation}
and
\begin{equation}\label{eq:dual-destra}
0 \to K \to \Oc_{X_r}^{k-1} \to \omega_C(5-r-c_1) \to 0
\end{equation}
we get that $h^1 \E^{\vee}(n)=0$ $\forall n \in \ZZ$ if and only if the map
$r_n:H^0\Oc_{X_r}^{k-1}(n) \to H^0\omega_C(5-r-c_1+n)$
is surjective $\forall n \in \ZZ$. The surjectivity of $r_n$ for
all $n<0$ is equivalent to $h^0\omega_C(4-r-c_1)=0$, for which a
sufficient condition is
$2g(C)-2<(r+c_1-4)\deg C$.
For $n=0$, the map $r_0$ is the one coming from the choice of
$k-1$ sections of $\omega_C(5-r-c_1+n)$, so it is surjective if
and only if we take a system of generators of
$H^0\omega_C(5-r-c_1+n)$. For $n>0$, observe that $r_n$ factors
through
\begin{equation*}
H^0 \Oc_{X_r}^{k-1} \otimes H^0 \Oc_{X_r}(n)
\to H^0 \omega_C(5-r-c_1) \otimes H^0 \Oc_C(n) \to 
H^0\omega_C(5-r-c_1+n).
\end{equation*}
If $C$ is projectively normal and $r_0$ is surjective, the
first map is surjective, and hence the surjectivity of $r_n$
becomes equivalent to the surjectivity of the second map,
which is precisely $\alpha_n$.

Finally, observe that the map $r_0$ is not injective (i.e.
condition (a) does not hold) is equivalent, by
(\ref{eq:dual-sinistra}) and (\ref{eq:dual-destra}) and the fact
that $c_1>0$, to the existence of a section of $\E^{\vee}$ mapping
to a nonzero section of $\Oc_{X_r}^{k-1}$. This is equivalent to
say that $\E$ has a trivial summand, which means that condition
(iii) in Definition \ref{def:star} dos not hold. On the other
hand, condition (b) is clearly equivalent, by (\ref{eq:basic}), to
$h^0\E(-1)=0$ and $h^0\E\ge k$, i.e. condition (i) in Definition
\ref{def:star}.
\end{proof}

We give now some examples of the above construction
applied to $X=X_4$ (for the rest of this section, unless
otherwise specified, $X$ will stand for $X_4$). We start
by giving a method to construct curves contained in a
hyperplane, which we will thus allow us to construct vector
bundles with $c_1=1$.

\begin{lemma}\label{lemma:degenere} Let ${\cal C}$ be a
family of degenerate curves in $\Pj^4$ such that some curve
of ${\cal C}$ is contained in a smooth degenerate (in $\Pj^4$) quartic
surface and it is not the complete intersection of this
surface and a hypersurface. Then the general quartic
hypersurface $X\subset\Pj^4$ contains a curve of ${\cal
C}$.
\end{lemma}

\begin{proof} Consider the set $\Sigma$ of degenerate
quartic surfaces of $\Pj^4$ containing a curve of ${\cal
C}$. We define the natural map
\begin{equation*}
\varphi:\Sigma\to ({\Pj^4})^*
\end{equation*}
associating to each quartic surface the unique hyperplane
containing it. Fix now a hyperplane $H$ of $\Pj^4$.
It is a standard fact in the Noether-Lefschetz theory
(see for example \cite{CGGH}), that the set of quartic
surfaces in $H$ containing a curve in ${\cal C}$ has
codimension one in the set of quartic surfaces of
$H$. In other words, the set $\varphi^{-1}(H)$ has dimension
$33$, and hence $\Sigma$ has dimension $37$.

Consider now the set $I$ of pairs $(S,X)$ where $S\in\Sigma$ and
$X$ is a quartic hypersurface of $\Pj^4$ containing $S$. Since the
fibers of the projection $I\to\Sigma$ are projective spaces of
dimension $35$, it follows that $I$ has dimension $72$.

We finally consider the second projection $p_2:I\to\Pj^{69}$
(where we identify $\Pj^{69}$ with the set of quartic
hypersurfaces of $\Pj^4$). The lemma will be proved if we
show that $p_2$ is surjective. Observe that, since any
smooth degenerate quartic surface in $\Pj^4$ is the
hyperplane section of some smooth quartic hypersurface of
$\Pj^4$, it follows from our hypotheses that there is a
smooth quartic hypersurface $X$ in the image of $p_2$.
Since a general hyperplane section of $X$ has its Picard
group generated by the hyperplane divisor (see for
example \cite{moi}), such a hyperplane section cannot
contain a curve of ${\cal C}$, and hence the set
$p_2^{-1}(X)$ has dimension at most three. Hence the
fiber of any element of the image of $p_2$ has
necessarily dimension three, which shows that $p_2$ is
surjective, as wanted.
\end{proof}

In order to apply the previous lemma we show the following:

\begin{lemma} \label{lem:key}
A smooth projectively normal space curve $C$ of degree $d$
and genus $g$ with
$d \geq g-1$
is contained in a smooth quartic surface in $\Pj^3$.
\end{lemma}

\begin{proof} The statement follows readily from the
more general results of \cite{mori}. We give however a direct
proof valid for our cases using the following simple standard
argument. We first observe, from Castelnuovo-Mumford's
criterion, that $\Ii_C(4)$ is globally generated. Indeed the
projective normality of $C$ provides the vanishing of $h^1
\Ii_C(3)$, while the vanishings of
$h^3\Ii_C(1)$ and $h^2 \Ii_C(2)$ come from the equalities
\begin{equation*}
h^3 \Ii_C(3)=h^2\Oc_C(3)=0,
\end{equation*}
\begin{equation*}
h^2 \Ii_C(2)=h^1 \Oc_C(2)=h^0 \omega_C(-2)=0
\end{equation*}
(the latter coming from the assumption $d \geq g-1$). 

Hence the linear system $|H^0 \Ii_C (4)|$  has no
base-points outside $C$ and therefore, by Bertini's theorem,
a general element of it is smooth outside $C$. By a well
known argument (see e.g. \cite{News}), from the exact
sequence
\begin{equation}\label{esatta-normale}
0 \to \Ii_C^2(4) \to \Ii_C(4) \to N_C^{\vee}(4) \to 0
\end{equation}
a surface in the linear system $|H^0 \Ii_C (4)|$ provides
a section of $N_C^{\vee}(4)$, and the singular points of the
surface belonging to $C$ are those in the zero locus of the
section. The sequence (\ref{esatta-normale}) proves that the
rank two vector bundle
$N_C^{\vee}(4)$ is generated by the global sections coming
from sections of $\Ii_C (4)$. Hence, a general such
section of $N_C^{\vee}(4)$ will be nowhere vanishing, which
implies that a general element of $|H^0 \Ii_C (4)|$ is
smooth also in the points of $C$. Such a general element
gives thus the wanted smooth quartic surface containing $C$.
\end{proof}

As a first application, we give a negative answer to
Question \ref{q:1}.

\begin{proposition} \label{ex:6-3}
There exists an ACM bundle on $X$ satisfying condition $\star$
with invariants $(4;1,6,4)$.
Moreover a general such
bundle is not extension class of any rank two ACM bundles.
\end{proposition}

\begin{proof}
Let $C$ be a smooth projectively normal curve
(hence non hyperelliptic) of degree $\deg C=6$ and genus $g(C)=3$ in $\Pj^3$.
By Lemma \ref{lem:key}, $C$ is contained in a smooth quartic
surface and by Lemma \ref{lemma:degenere} a general quartic
hypersurface in $\Pj^4$ contains a degenerate curve $C$ of
degree $6$ and genus $3$. Since $h^0 \omega_C=3$,
by Proposition \ref{prop:vogelaar}
$C$ defines a rank four vector bundle given
by
\begin{equation} \label{eq:f2}
0 \to \Oc_X^3 \to \E \to \Ii_C(1) \to 0
\end{equation}
with $c_1(\E)=1$, $c_2(\E)=6$, $c_3=4$ (by
(\ref{eq:genus})) and
$h^0\E=4$. Proposition \ref{pro:criterio-ACM} applies once
we show the map
\begin{equation*}
\alpha_1:H^0 \omega_C \otimes H^0 \Oc_C(1) \to H^0 \omega_C(1)
\end{equation*}
is surjective, which follows
by Castelnuovo's Lemma (see
\cite{ACGH} pg.151 or \cite{AS} theorem (1.6)). Hence $\E$
is ACM and satisfies condition $\star$. 

Finally any such bundle is not an
extension class of any two rank two ACM bundles on $X$.
Suppose to the contrary that there exists
a non trivial extension class
\[
0 \to \E' \to \E \to \E'' \to 0
\]
with $\E'$ and $\E''$ ACM of rank two. Since $\E$ satisfies
condition $\star$, then both $\E'$ and $\E''$ are normalized, i.e.
$b(\E')=b(\E'')=0$ and $h^0 \E' \cdot h^0 \E'' \geq 1$. Then a direct
computation shows that equations (\ref{eq:wit}) have not integral
solutions for Chern classes of $\E'$ and $\E''$ as in Proposition
\ref{prop:list} and we are done.
\end{proof}

\begin{remark} \normalfont{
Alternatively, to show the existence of space curves of degree $6$
and genus $3$ on general $X$ one can start with a curve of degree
$10$ as in Example \ref{ex:4,2} and then taking the residual curve
to it in a complete intersection $(2,2,4)$.}
\end{remark}

Similar to the previous case, we have also the following:

\begin{proposition} \label{ex:3-i}
There exists an ACM bundle on $X$ satisfying condition $\star$
with invariants $(3;1,5,2)$.
\end{proposition}

\begin{proof}
Let $C$ be a smooth curve of type $(2,3)$ in a smooth quadric surface
$C \subset Q \subset \Pj^3$. By Lemma \ref{lem:key}, there
exists a smooth quartic surface containing $C$. Hence, by
Lemma \ref{lemma:degenere}, a general quartic hypersurface
contains a curve $C$ of degree $5$ genus $2$.
Since $h^0 \omega_C=2$,
by Proposition \ref{prop:vogelaar} we get a rank three vector
bundle
$\E$ fitting in the exact sequence
\[
0 \to \Oc_X^2 \to \E \to \Ii_C(1) \to 0
\]
with $c_1=1$, $c_2=5$, $c_3=2$ and $h^0 \E=3$. 
Using Castelnuovo's Lemma 
it is easy to
check that $C$ satisfies the conditions of Proposition
\ref{pro:criterio-ACM} 
and
hence $\E$ is ACM and satisfies condition $\star$.

\end{proof}

\begin{remark}
\normalfont{
A similar liaison argument as above shows the existence
of rational space cubic curves in $X$. Specifically, let $Q$
be the unique quadric surface containing a quintic $C\subset
X$ of genus $2$ as in Proposition \ref{ex:3-i}. Then
$C$ is residual to a rational cubic curve $D \subset Q \cap X$.
}
\end{remark}

\section{Bounds of the Chern classes} \label{sec:bounds}

In this section we will prove several restrictions for the
Chern classes of ACM vector bundles, with special
attention to the case $r=4$. This will yield the list of
possible Chern classes stated in Theorem
\ref{thm:primoelenco} for the cases $k=3,4$.
We start with bounds for $c_1$ and $c_2$. The proof of the
bound for $c_1$ is a straightforward extension of the method
used in \cite{Ma1} for the rank two case.

\begin{lemma} \label{lem:bound-generale}
Let $(\E,C)$ be an ACM pair on $X_r$ satisfying
condition $\star$ and having invariants $(k;c_1,c_2,c_3)$.
Then
\begin{equation}\label{primobound:c1}
1\le c_1\le \frac{k(r-1)}2
\end{equation}
and
\[
c_2\le \frac{r}2 c_1^2-\frac{r(r-2)}{2} c_1+\frac{r(r-1)(r-2)}{6} k.
\]
\end{lemma}

\begin{proof}
By the exact sequence (\ref{eq:basic}), since $h^0\E \geq k$
it follows $h^0\Ii_C(c_1) \geq 1$ and hence $c_1\geq 1$.
Let $H$ be a general hyperplane section of $X_r$.
Taking cohomology in the exact sequence
\[
0 \to \E(-2) \to \E(-1) \to \E_H(-1) \to 0
\]
(and its twists by any $\Oc_{X_r}(l)$), we get that
$\E_H$ is an ACM bundle on $H$. Since  $h^0\E(-1)$ (because
$\E$ normalized) we also get $h^0\E_H(-1)=0$. Therefore,
\[
0\le h^2\E_H(-1)  =\chi\E_H(-1) =-c_2 +
\frac{r}{2}c_1^2-\frac{r(r-2)}{2}c_1 
+\frac{r(r-1)(r-2)}{6}k
\]
and the bound for $c_2$ follows.
Similarly, if $Y$ is a general hyperplane section of
$H$, we have an exact sequence
\[
0 \to \E_H(-2) \to \E_H(-1) \to \E_Y(-1) \to 0
\]
which implies $h^0\E_Y(-1)=0$. Therefore
$0 \ge\chi\E_Y(-1)=r(c_1-k(r-1)/2)$
and the upper
bound for $c_1$ follows immediately.
\end{proof}

In the case of quartic hypersurfaces $X_4$, which is the
case we are interested in, we have stronger restrictions:

\begin{proposition}\label{pro:k-1}
Let $(\E,C)$ be an ACM pair satisfying condition
$\star$ on $X_4$ and having invariants $(k;c_1,c_2,c_3)$.
Then
\begin{equation}\label{eq:c3}
c_3=-\frac{4}{3}c_1^3+2c_1^2-\frac{14}{3}c_1+c_1c_2-c_2+2k
\end{equation}
and
\begin{equation}\label{eq:genus12}
g(C)=-\frac{2}{3}c_1^3+c_1^2-\frac{7}{3}c_1+1+(c_1-1)c_2+k.
\end{equation}
Moreover
\begin{equation}\label{bound:c1}
1 \leq c_1 \leq \frac{3k}2.
\end{equation}
\begin{equation}\label{eq:bound-c2}
 2c_1^2-2c_1+k\le c_2\le\min\{2c_1^2-4c_1+4k,2c_1^2+k\}
\end{equation}
and if $c_1>1$ then
$c_2 \geq 2c_1^2-4c_1+8$.
\end{proposition}

\begin{proof}
By hypothesis $h^0\E(-1)=0$, and also $h^0\E^\vee=0$ because
$\E$ has not trivial summands (see the proof of Proposition
\ref{pro:criterio-ACM}). Since $h^3\E(-1)=h^0\E^\vee$ (this
is the point in which we use $r=4$) it follows
$\chi\E(-1)=0$, which provides (\ref{eq:c3}) by using
Riemann-Roch. Now (\ref{eq:genus12}) comes from
(\ref{eq:genus}) performing the substitution of $c_3$ given
by (\ref{eq:c3}).

The bound (\ref{bound:c1}) and first upper bound for
$c_2$ in (\ref{eq:bound-c2}) are the ones of Lemma
\ref{lem:bound-generale} for $r=4$. For the other upper
bound for $c_2$, since $h^3\E=h^0\E^\vee(-1)=0$ (again
because $r=4$) and $h^0\E\ge k$ by hypothesis, if follows
that
$\chi(\E)\ge k$. By Riemann-Roch, using the substitution
(\ref{eq:c3}) we get
$\chi(\E)=-c_2+2c_1^2+2k$
which yields the wanted upper bound
for $c_2$.

For the lower bound of $c_2$, we take a general
linear projection
$\pi:X\to\Pj^3$ from a point of $\Pj^4$. It is immediate to
observe that $\pi_*\E$ is an ACM vector bundle of rank $4k$,
which should split completely by Horrocks theorem. Since
$h^0\E(-1)=0$, $\pi_*\E$ must contain a direct summand
$\Oc_{X_4}^{h^0\E}$, and similarly a direct summand
$\Oc_{X_4}(-2)^{h^3\E(-2)}$, since $h^3\E(-1)=0$. Therefore,
$h^0\E+h^3\E(-2)\le 4k$, which gives the wanted inequality
when using Riemann-Roch and the substitution
\eqref{eq:c3}.

Finally, let us prove the last lower bound for $c_2$. Since
$h^0
\Ii_C(c_1-2)=0$ (by Proposition
\ref{pro:criterio-ACM}(b)), it follows that
\[
h^0\Oc_C(c_1-2)=h^0\Oc_{X_4}(c_1-2)
=\frac{2}{3}c_1^3-3c_1^2+\frac{19}{3}c_1-5.
\]
Since by Riemann-Roch we have
\[
\chi\Oc_C(c_1-2)=\frac{2}{3}c_1^3-c_1^2+\frac{7}{3}c_1-k-c_2
\]
we derive
$h^0\omega_C(2-c_1)=c_2-2c_1^2+4c_1-5+k$.
The proof concludes by applying the bilinear lemma (see \cite{cl}) to the
nondegenerate bilinear map
\[
H^0\omega_C(1-c_1)\otimes H^0\Oc_C(1)\to H^0\omega_C(2-c_1)
\]
using that $h^0\omega_C(1-c_1)=k-1$ (by Proposition
\ref{pro:criterio-ACM}(a)) and that $h^0\Oc_C(1)\ge 5$ if $c_1>1$,
because $C$ is not contained in any hyperplane (by Proposition
\ref{pro:criterio-ACM}(b)).
\end{proof}

We will study separately the cases $c_1=1$ and
$2 \leq c_1 \leq 3k/2$. Proposition \ref{pro:k-1} gives
immediate results for vector bundles of rank $k \geq 3$ with
$c_1=1$:

\begin{proposition}\label{pro:c1=1}
Let $(\E,C)$ be a pair on $X_4$ satisfying condition $\star$,
where $\E$ is ACM of rank $k\ge3$ and $c_1=1$. Then
$c_2=k+2$, $c_3=2k-4$ and
$C$ is a curve of degree $k+2$ and genus $g(C)=k-1$.
\end{proposition}

\begin{proof}
We immediately get from Proposition \ref{pro:k-1} $g(C)=k-1$
and $\deg C=c_2\le k+2$. Hence, $C$ cannot be a plane curve
and therefore $h^0\Ii_C(1)=1$ (it is at least $1$ by
Proposition
\ref{pro:criterio-ACM}). This implies ($h^3\E=0$, as proved
in Proposition \ref{pro:k-1})
$\chi(\E)=h^0\E=k-1+h^0\Ii_C(1)=k$. Since by Riemann-Roch
$\chi(\E)=-c_2+2k+2$
it follows $c_2=k+2$, which is the degree of $C$.
Finally, from (\ref{eq:c3}) we get $c_3=2k-4$.
\end{proof}

For the rank three case when $c_1>1$ we have the following
result if $c_1=3,4$:

\begin{proposition} \label{prop:bound-migliore}
Let $(\E,C)$ be a pair satisfying condition $\star$, where
$\E$ is ACM of rank $k=3$. If $c_1\ge 3$ then
$2c_1^2-4c_1+11\le c_2\le 2c_1^2-4c_1+12$.
\end{proposition}

\begin{proof}
The upper bound follows from Proposition \ref{pro:k-1}. For
the lower bound, observe that $h^3 \E^\vee(1)=h^0 \E(-2)=0$,
and thus
$h^0 \E^\vee(1)=\chi \E^\vee(1)=2c_1^2-4c_1-c_2+12$.
Hence it
will be enough to prove $h^0 \E^{\vee}(1)\leq 1$.

Assume for contradiction $h^0 \E^{\vee}(1) \geq 2$ and let
$s_1,s_2$ be two independent global sections. 
Then $\wedge^2 (\E^{\vee}(1)) \cong
\E(2-c_1)$ contains the nonzero section defined by $s_1
\wedge s_2$, which is absurd since $\E$ is normalized and
$2-c_1<0$. \end{proof}

\begin{remark}\normalfont{
The above Proposition leaves only two possibilities for
$c_2$, depending on whether $h^0 \E^{\vee}(1)$ is $0$ or
$1$. When $h^0 \E^{\vee}(1)=0$ (i.e. $c_2=2c_1^2-4c_1+12$),
it follows that both $\E(-1)$ and its dual have no global
sections, so it follows (as in \cite{OSS} Remark 1.2.6) that
$\E(-1)$, and hence $\E$, is stable.
Observe that if $c_1=3$ we get $c_2=18$, and if $c_1=4$ then
$c_2=28$.
}\end{remark}

Propositions \ref{pro:k-1}, \ref{pro:c1=1} and
\ref{prop:bound-migliore} give immediately the following
possibilities for the Chern classes of ACM vector bundles of
rank three or four over $X_4$:

\begin{proposition}
Let $(\E,C)$ be an ACM pair satisfying condition $\star$
where $\E$ is of rank $k=3$. Then the possibilities
for the Chern classes of $\E$ and the genus $g=g(C)$ are:
\begin{enumerate}
\item
$c_1=1$, $c_2=5$, $c_3=2$, $g(C)=2$;
\item
$c_1=2$, $8 \leq c_2 \leq 11$, $c_3=c_2-6$,
$g(C)=c_2-2$;
\item
$c_1=3$, $17 \leq c_2 \leq 18$,  $c_3=2c_2-26$,
$g(C)=2c_2-12$;
\item
$c_1=4$, $27 \leq c_2 \leq 28$, $c_3=3c_2-66$,
$g(C)=3c_2-32$.
\end{enumerate}
\end{proposition}

\begin{proposition}
Let $(\E,C)$ be an ACM pair satisfying condition $\star$
where $\E$ is of rank $k=4$. Then the possibilities for the
Chern classes of $\E$ and the genus $g=g(C)$ are:
\begin{enumerate}
\item
$c_1=1$, $c_2=6$, $c_3=4$, $g(C)=3$;
\item
$c_1=2$, $8 \leq c_2 \leq 12$, $c_3=c_2-4$,
$g(C)=c_2-1$;
\item
$c_1=3$, $16 \leq c_2 \leq 22$, $c_3=2c_2-24$,
$g(C)=2c_2-11$;
\item
$c_1=4$, $28 \leq c_2 \leq 32$, $c_3=3c_2-64$,
$g(C)=3c_2-31$;
\item
$c_1=5$, $44 \leq c_2 \leq 46$, $c_3=4c_2-132$,
$g(C)=4c_2-65$;
\item
$c_1=6$, $c_2=64$, $c_3=84$,
$g(C)=5c_2-117$.
\end{enumerate}
\end{proposition}

We end this section with some remarks and comments.
We notice the following result can be
derived using similar arguments as in the above cases of rank three:

\begin{corollary}
There are no rank three ACM vector bundles satisfying condition
$\star$ on a general quartic hypersurface in $\Pj^n$, $n \geq 5$.
In other words, rank three ACM bundles satisfying condition $\star$
on a general quartic threefold $X_4 \subset \Pj^4$ do not extend
to a general quartic hypersurface $W \subset \Pj^5$ having $X_4$
as a hyperplane section.
\end{corollary}

Since in \cite{CM3} we showed that all rank two ACM bundles on a
general quartic $W \subset \Pj^n$, $n \geq 5$, splits, then a
natural question, which we hope to consider later, is the
following:

\begin{question}
Which is the minimum rank for a positive dimensional family of ACM
bundles satisfying condition $\star$ on a general quartic fourfold
$W \subset \Pj^5$.
\end{question}

\

\

\noindent Enrique Arrondo \par\noindent
Departamento de \'Algebra,
Facultad de Ciencias Matem\'aticas,
Universidad Complutense de Madrid,
28040 Madrid (Spain) \par\noindent
email: Enrique\_Arrondo@mat.ucm.es

\par\bigskip

\noindent Carlo G. Madonna 
\par\noindent Departamento de Matem\'aticas, Instituto
de Matem\'aticas y F\'{\i}sica Fundamental, CSIC, C/ Serrano
121, 28006 Madrid (Spain)
\par\noindent 
{\it current address}: 
\par\noindent
Departamento de \'Algebra, Facultad de Ciencias Matem\'aticas 
Universidad Complutense de Madrid, 28040 Madrid (Spain)
\par\noindent email: carlo.madonna@mat.csic.es


\begin{thebibliography}{ADSE}

\bibitem{ACGH} E. Arbarello, M. Cornalba,
P.A. Griffiths and J. Harris, Geometry of
Algebraic Curves, Springer, 1985.

\bibitem{AS} E. Arbarello and E. Sernesi,
\emph{ Petri's approach to the study of the
ideal associated to a special divisor},
Inventiones Math. \textbf{49}(1978), 99--119.

\bibitem{arrondo} E. Arrondo, \emph{A home-made Hartshorne-Serre
correspondence}, Rev. Mat. Complut.  {\bf 20} (2007),  no. 2, 423--443.

\bibitem{AC} E. Arrondo and L. Costa,
\emph{Vector bundles on Fano 3--folds without
intermediate cohomology}, Comm. Algebra {\bf
28} (2000), no. 8, 3899--3911.

\bibitem{AF} E. Arrondo and D. Faenzi,
\emph{Vector bundles with no intermediate
cohomology on Fano threefolds of type
$V_{22}$}, Pacific J. Math.  \textbf{225}  (2006),  no. 2,
201--220.

\bibitem{CGGH} J. Carlson, M. Green, P.
Griffiths and J. Harris, \emph{Infinitesimal
variations of Hodge structures (I)}, Comp.
Math. \textbf{50} (1983), 105--205.

\bibitem{cm1} L. Chiantini and C. Madonna,
\emph{ACM bundles on a general quintic
threefold}, Matematiche (Catania)
\textbf{55}(2000), no.2, 239--258.

\bibitem{CM2} L. Chiantini and C. Madonna,
\emph{A splitting criterion for rank 2 vector bundles on a general
sextic threefold}, Internat. J. Math.
\textbf{15}(2004), no.4, 341--359.

\bibitem{CM3} L. Chiantini and C.Madonna,
\emph{
ACM bundles on general hypersurfaces in $\Pj^5$ of low degree},
Collect. Math. Vol. \textbf{56}(2005), no. 1, 85-96.

\bibitem{EH} D. Eisenbud, Commutative algebra
with a view toward algebraic geometry,
Springer 1999.

\bibitem{cl} D. Eisenbud, J. Koh, and
M. Stillman, \emph{ Determinantal Equations for
Curves of High Degree} Amer. J. Math. Vol.
\textbf{110}(1988), no. 3, pp. 513-539.

\bibitem{F} D. Faenzi,
\emph{Bundles over the Fano threefold $V_5$},  Comm. Algebra  \textbf{33}(2005), no. 9, 3061--3080.

\bibitem{fulton} W. Fulton, Intersection
theory, Springer 1998.

\bibitem{hrs} J. Harris, M. Roth and J. Starr,
\emph{Curves of small degree on cubic
threefolds}, Rocky Mountain J. Math.
\textbf{35}(2005), no. 3, 761--817.

\bibitem{Hor} G. Horrocks, \emph{Vector
bundles on the punctured spectrum of a local
ring}, Proc. London Math. Soc. \textbf{14}
(1964), 689-713.

\bibitem{Ili} A. Iliev and D. Markushevich,
\emph{Quartic 3-folds: pfaffians, vector
bundles, and half-canonical curves}, Mich.
Math. J. \textbf{47} (2000), 385--394.

\bibitem{Ili-Man} A. Iliev and L. Manivel,
\emph{Pfaffian lines and vector bundles on
Fano threefolds of genus $8$}, J. Alg. Geom.,
to appear.

\bibitem{KRR} N.M. Kumar, A.P. Rao and G.V. Ravindra, \emph{
Arithmetically Cohen-Macaulay bundles on three dimensional hypersurfaces},
Int. Math. Res. Not. IMRN  2007,  no. 8, Art. ID rnm025, 11 pp.

\bibitem{Ma1} C. Madonna, \emph{A splitting
criterion for rank 2 vector bundles on
hypersurfaces in $\Pq$}, Rend. Sem. Mat. Univ.
Pol. Torino {\textbf 56} (1998), no.2, 43--54.

\bibitem{Ma2} C. Madonna, \emph{Rank--two
vector bundles on general quartic
hypersurfaces in $\Pq$}, Rev. Mat. Complut.
{\bf XIII} (2000), num.2, 287--301.

\bibitem{roum} C.G. Madonna, \emph{ACM vector
bundles on prime Fano threefolds and complete
intersection Calabi Yau threefolds}, Rev.
Roumaine Math. Pures Appl. \textbf{47}(2002)
no.2, 211-222.

\bibitem{MadCEJM} C. Madonna, \emph{Rank 4
vector bundles on the quintic threefold}, CEJM
\textbf{3} (2005), no.3, 404--411.

\bibitem{mori} S. Mori, \emph{On degrees and genera of curves on smooth quartic surfaces in $\Pj^3$}, 
Nagoya Math. J. \textbf{96}(1984), 127--132. 

\bibitem{moi} B.G. Mo\u\i \v sezon, \emph{
Algebraic homology classes on algebraic
varieties}, Izv. Akad. Nauk SSSR
Ser. Mat. \textbf{31} (1967) 225--268.

\bibitem{News} P.E. Newstead, \emph{A space curve whose normal bundle is stable},
J.London Math. Soc. \textbf{28} (1983), no.2, 428--434.


\bibitem{OSS} C. Okonek, M. Schneider and
H. Splinder, Vector bundles on complex
projective spaces, Birkhauser, 1980, Boston.

\bibitem{Ott} G. Ottaviani, \emph{Some
extension of Horrocks criterion to vector
bundles on Grassmannians and quadrics}, Ann.
Mat. Pura Appl. (4)
\textbf{155} (1989), 317--341.

\bibitem{vog} A. Vogeelar, \emph{Constructing
vector bundles from codimension-two
subvarieties}, PhD thesis, Leiden, 1978.

\end{thebibliography}
\end{document}